\documentclass[11pt]{article}
\usepackage{amsfonts}
\usepackage{mathrsfs}
\usepackage{amssymb}
\usepackage{amsmath}
\usepackage[all]{xy}
\newtheorem{Lemma}{Lemma}
\newtheorem{Theorem}{Theorem}
\newtheorem{Corollary}{Corollary}
\newtheorem{Definition}{Definition}

\def\Coker{\mathop{\rm Coker}\nolimits}
\def\Tr{\mathop{\rm Tr}\nolimits}

\def\Ext{\mathop{\rm Ext}\nolimits}
\def\Hom{\mathop{\rm Hom}\nolimits}
\def\Tor{\mathop{\rm Tor}\nolimits}

\title{\Large \bf
Local algebras with radical cubic zero are  PCM-free
}
\author{Luo Rong\thanks{\it E-mail:luorong@home.swjtu.edu.cn}$\quad\quad$ \\{\small \it College  of Mathematics, Southwest
Jiaotong University, Chengdu 610031, P. R. China}}
\usepackage{amssymb}
\date{}
\begin{document}
\baselineskip=18pt \maketitle

\begin{abstract}
An artin algebra is said to be PCM-free if every finitely generated Gorenstein projective module with a projective submodule  is projective. In this paper, we show  that  artin local algebras with radical cubic zero are PCM-free.
\end{abstract}

\section{Introduction}

~~~~~Throughout this paper, $R$ is a commutative artin ring,  $\Lambda$ is an  artin algebra, that is an $R$-algebra which is a
finitely generated $R$-module and mod$\Lambda$ the category of finitely generated left $\Lambda$-modules. We recall from [AM] a complex ${\bf P^*}$ of projective
$\Lambda$-modules is totally acyclic if it is acyclic and for each projective $\Lambda$-module $Q$
the $\Hom$ complex $\Hom_\Lambda({\bf P^*,Q})$ is acyclic. A $\Lambda$-module $M$ is called a (finitely generated)
Gorenstein projective module ([EJ]) provided that there is a totally acyclic
complex ${\bf P^*}$ such that the zeroth cocycle $Z^0{\bf P^*}$ is isomorphic to $M$. In this case,
the complex ${\bf P^*}$is called a complete resolution of $M$. In the literature, Gorenstein
projective modules are also called modules of G-dimension zero ([AuB]).

We denote $G(\Lambda)$ as the full subcategory of mod$\Lambda$
consisting of all Gorenstein projective modules, and  the full
subcategory of mod$\Lambda$ being composed of all projective module
by $P(\Lambda)$. It is well known that a projective module is
Gorenstein projective. That is, $P(\Lambda)\subset G(\Lambda)$. An
artin algebra is called CM-free provided that
G$(\Lambda)=P(\Lambda)$, that is, its all finitely generated
Gorenstein projective modules are projective. An extreme case has
been considered by Xiaowu Chen in [Ch]: a connected artin ring with
radical square zero is either CM-free or self-injective.

In this paper, another extreme case is considered. We recall from
[T] that a Gorenstein projective module is said to be P-Gorenstein
projective if it contains a nonzero-projective submodule. An artin
algebra is PCM-free if every finitely generated P-Gorenstein
projective module is projective. Note that a CM-free  algebra is
PCM-free. However, the converse is not true in general.

For an algebra $\Lambda$, denote by $J$ its Jacobson radical. The
algebra $\Lambda$ is said to be with radical cubic zero provided
that $J^3=0$. Let us remark that there exists an upper triangular
matrix algebra which is not CM-free is such a Gorenstein algebra
that its Jacobson radical satisfies the radical cubic zero and
square nonzero; see the remarks after [Ch]. Here we would like to
address a problem: for an algebra with the radical cubic zero, when
are  finitely generated Gorenstein projective modules projective?

 The aim of the paper is to show that for an artin local
algebra with radical cubic zero the study of its finitely generated
P-Gorenstein projective modules always belongs to the one extreme
case as follows.

\begin{Theorem}\label{Th1}
Let $(\Lambda, S)$ be an artin local algebra with radical cubic zero. Then $\Lambda$ is  PCM-free.
\end{Theorem}

We draw an immediate consequence of Theorem 1. Recall that a local algebra
$\Lambda$ is $n$-Gorenstein provided that  G-dimension of the simple module is less than or equal to $n$.
A local algebra $\Lambda$ is $n$-regular if projective dimension of the simple module is less than or equal to $n$(see [C]).

\begin{Corollary}
Let $(\Lambda, S)$ be an artin local algebra with radical cubic zero. If $\Lambda$ is $1$-Gorenstein, then it is either $1$-regular or self-injective.
\end{Corollary}
{\bf Proof.} It suffices to notice the following fact: If G-dimension of $S$ is less than or equal to $1$, then $S$ admits a surjective Gorenstein projective precover
$\psi: G\twoheadrightarrow S$ where
$K=Ker\psi$ satisfies $K$ is projective and $G$ is indecomposable. That is, there exists an exact sequence $0\to K\to G\to S\to 0$. If $K\neq 0$, by the theorem 1, we know that
$G$ is projective. So projective dimension of $S$ is less than or equal to $1$. That is, $\Lambda$ is $1$-regular. If $K=0$, $S$ is Gorenstein projective. So,
$\Lambda$ is self-injective.
\vspace{0.2cm}

In the next section, we start by recalling the definition of irreducible morphisms,
give several preliminary involving  properties of Gorenstein projective modules, and
finally prove the above theorem.

\section{ Proof of the Theorem}

~~~~~In this section we present the proofs of Theorem $1\label{Thii}$.
 Let $\Lambda$ be an artin algebra. Recall that for each $\Lambda$-module $M$, its syzygy module
$\Omega(M)$ is defined to be the kernel of its projective cover $P\to M$. Recall that in
a short exact sequence $0\to M_1\to P\stackrel{p}\rightarrow M\to 0$ with $P$ projective, we have
$M_1\cong \Omega(M)\oplus Q$ for a projective module $Q$; moreover, $Q=0$ if and only if $p$ is a
projective cover. Here we state several properties of Gorenstein projective modules for later use. For the proofs,
we refer to [Ch, Lemma 2.1, 2.2].
\begin{Lemma}\label{L1}
Let $M$ be a Gorenstein projective $\Lambda$-module which is indecomposable
and non-projective. Then $\Omega(M)$ is also an indecomposable non-projective Gorenstein projective $\Lambda$-module.
\end{Lemma}

Recall that $J$ denotes the Jacobson
radical of $\Lambda$.
\begin{Lemma}\label{L2}
 Let $M$ be a Gorenstein projective $\Lambda$-module without projective direct
summands. Assume that $J^n=0$ for $n\geq 2$. Then $J^{n-1}M=0$ and $J^{n-2}M$ is semi-simple.
\end{Lemma}

Now we introduce the notion of an irreducible morphism and state its  property for later use. For the proof of its property,
we refer to [AuR1, Proposition 4.1].
\begin{Definition}
  A morphism $f: B\to C$ in
the category {\rm mod}$\Lambda$  is said to be Irreducible if $(a)$ $f$
 is not a splittable proper epimorphism or a splittable proper
monomorphism;  $(b)$  given a commutative diagram
$$\xymatrix{&Z\ar[dr]^{h}\\
B\ar[rr]^f\ar[ur]^g&&C}$$  either $h$  is splittable
epimorphic or $g$  is splittable  monomorphic.
\end{Definition}
\begin{Lemma}\label{L3}
Let $M$ be an indecomposable $G$-projective module. The following conditions
are equivalent:

$(1)$ There exists an irreducible monomorphism $M\to Q$ for a projective $Q$.

$(2)$ $M$ is isomorphic to a summand of $JQ$ for a projective $Q$.
\end{Lemma}

Let $P_1\stackrel{f}\longrightarrow P_0\to M\to 0$ be a minimal projective resolution of a module $M$ in mod$\Lambda$. We call $\Coker f^*$
the Auslander transpose of $M$, and denote it by $\Tr M$([AuR]).  Consider the stable category $\underline{{\rm mod}}\Lambda$ of mod$\Lambda$, by the [AuR, Proposition 2.2], we have the  following lemma.
\begin{Lemma}\label{L4}
Let $M, N$ be in {\rm mod}$\Lambda$. Then there exists an isomorphism $\Tor_1^\Lambda(\Tr M,N)\cong \underline{\Hom}_\Lambda(M,N)$.
\end{Lemma}

Let $D(-)$ be the nature functor $\Hom_R(-, E(R/r))$ where $r$ is Jacobson
radical of $R$ and $E(R/r)$ is the injective envelope.  We say that an exact sequence $0\longrightarrow
A\stackrel{f}\longrightarrow B \stackrel{g}\longrightarrow C\longrightarrow 0$ in mod$\Lambda$ is almost split if $f$ is left almost split
and $g$ is right almost split (see [AuR]). It is well known that there  exists a unique almost split sequence for a non-injective indecomposable $\Lambda$-module $M$ and the form is $0\to M\to N\to \Tr DM\to 0$. Now the time has come when we can prove our theorem \ref{Th1}.
\vspace{0.2cm}

{\bf The proof of Theorem \ref{Th1}} Let $PG(\Lambda)$ be the subcategory of mod$\Lambda$ composed of all
finitely generated $P$-Gorenstein projective modules. Assume that
$\Lambda$ is not PCM-free. Take a $M$ in $PG(\Lambda)$ to be
indecomposable and non-projective. Set $Q$ to be the
nonzero-projective submodule $M$. Considering the monomorphism $0\to
Q\to M$, there is the monomorphism $0\to JQ\to JM$. Note that
$J^3=0$. By the lemma 2  we have $JM$ is semisimple and then $JQ$ is
semisimple. Let $S$ be the simple submodule of $JQ$. Then there
exists an indecomposable projective module $P$ such that the simple
module $S$ is isomorphic to the summand of $JP$. It follows
immediately from the lemma 3  that there exists an irreducible
morphism $f:S\to P$. Take a short exact sequence $0\to S\to E\to \Tr
DS\to 0$ such that it is the almost split sequence of $S$. Then
there exist the morphisms $h:E\to P$ and $g:\Tr DS\to P/S$ such that
the following diagram
$$\xymatrix{0\ar[r]&S\ar[r]\ar@{=}[d]&E\ar[r]\ar[d]^h&\Tr DS\ar[r]\ar[d]^g&0\\
0\ar[r]&S\ar[r]^f&P\ar[r]&P/S\ar[r]&0}$$ commutates. Note that $f$
is irreducible, $h$ is split. This gives rise to the almost split
exact sequence of $S$: $$0\to S\to P\oplus B\to \Tr DS\to 0$$ Hence
we have a commutative diagram

$$\xymatrix{&&0\ar[d]&0\ar[d]\\
&&B\ar@{=}[r]\ar[d]&B\ar[d]\\
0\ar[r]&S\ar@{=}[d]\ar[r]&P\oplus B\ar[r]\ar[d]&\Tr DS\ar[r]\ar[d]&0\\
0\ar[r]&S\ar[r]&P\ar[r]\ar[d]&P/S\ar[r]\ar[d]&0\\
&&0&0}$$
This gives us the commutative diagram
$$\xymatrix{0\ar[r]&S\ar[r]\ar[d]^\theta&P\ar[r]\ar[d]^\sigma&P/S\ar[r]\ar@{=}[d]&0\\
0\ar[r]&B\ar[r]&\Tr DS\ar[r]&P/S\ar[r]&0}$$
That $\sigma$ is irreducible follows from the relationship between irreducible morphisms and almost split sequence.

{\bf A):} If $\theta=0$, by the indecomposable module $\Tr DS$, $\Tr DS\cong P/S$. And we have the almost split sequence $$\eta: 0\to S\to P\to \Tr DS\to 0$$
If $M$ is isomorphic to $S$, by a projective submodule $Q$ of $M$, $M$ is projective. It is contradicted with our assumption. Suppose that $M$ is not isomorphic to $S$.
Note that the exact sequence $\eta$ is an almost
split sequence, we know that $\underline{\Hom}_\Lambda(S,M)=0$. It follows from Lemma \ref{L4} that
$$\xymatrix{\Tor_1^\Lambda(\Tr S, M)=0 & {\rm and}& \Ext_\Lambda^1(\Tr S, DM)\cong D\Tor^\Lambda_1(\Tr S,M)=0}$$
Let $0\to M\to P_1\to M_1\to 0$ be an exact sequence with a projective module $P_1$.
This means that $\Ext_\Lambda^2(\Tr S, DM_1)=0$. Use the functor $\Hom_\Lambda(-,\Lambda)$ to the almost split sequence $\eta$,
by the non-projective module $S$, we obtain the short exact sequence $$0\to DS\to P'\to \Tr S\to 0$$ with the projective module  $P'$.
 Since $\Ext_\Lambda^2(\Tr S, DM_1)=0$, we have that
\begin{eqnarray*}
  \Ext_\Lambda^1(M_1, S)\cong\Ext_\Lambda^1(DS, DM_1) \cong\Ext_\Lambda^2(\Tr S,DM_1)=0
\end{eqnarray*}
Therefore, $M_1$ is projective. It is contradicted with the non-projective module $M$.

{\bf B):} If $\theta\neq 0$, we have the monomorphism $\sigma: P\to \Tr DS$. Noting the almost split sequence $0\to S\to E\to \Tr DS\to 0$, we
have a commutative diagram
\begin{eqnarray*}
  \xymatrix{&P\ar[d]^\sigma\ar[dl]_\alpha\\E\ar[r]&\Tr DS\ar[r]&0.}
\end{eqnarray*}
Since $\sigma$ is irreducible, it follows that $\alpha$ is a splittable monomorphism. That is, there exists a epimorphism $\beta:E\to P$
such that $\beta\alpha=1_P$. Consider the commutative diagram
\begin{eqnarray*}
  \xymatrix{&&0\ar[d]&0\ar[d]\\&&P\ar@{=}[r]&P\ar[d]\\0\ar[r]&S\ar[r]^\gamma&E_1\oplus P\ar[u]^\beta\ar[r]&\Tr DS\ar[r]&0}
\end{eqnarray*}
with $E\cong (P\oplus E_1)$. If $\beta\gamma\neq 0$, it is contradicted with the monomorphism $\sigma$.   Hence we have the isomorphisms
$E/S\cong (E_1\oplus P)/S\cong (E_1/S)\oplus P\cong \Tr DS$.
This implies that $P=0$. The result contraries to the assumption of the theorem. This contradiction completes the
proof of the theorem.

\end{document}